\documentclass[10pt,leqno]{article}
\usepackage{amssymb,amsmath,amscd,dsfont}
\newcommand{\op}[1]{\ensuremath{\operatorname{#1}}}
\newcommand{\ol}[1]{\ensuremath{\overline{#1}}}
\newcommand{\R}{\ensuremath{\mathds{R}}}
\newcommand{\N}{\ensuremath{\mathds{N}}}
\newcommand{\id}{\ensuremath{\operatorname{id}}}
\newcommand{\pr}{\ensuremath{\operatorname{pr}}}
\newcommand{\supp}{\ensuremath{\operatorname{supp}}}
\newcommand{\se}{\ensuremath{\nobreak\subseteq\nobreak}}
\newcommand{\from}{\ensuremath{\nobreak:\nobreak}}
\renewcommand{\to}{\ensuremath{\nobreak\rightarrow\nobreak}}

\usepackage[a4paper,hscale=0.7,vscale=0.82,hmarginratio=1:1,includehead,headheight=14pt]{geometry}
\usepackage{fancyhdr}
\fancypagestyle{plain}{\fancyhead{}\fancyfoot{}}

\newcommand{\shortTitle}{}

\usepackage{amsthm}
\theoremstyle{definition}
\newtheorem{definition}{Definition}[section]
\newtheorem{remark}[definition]{Remark}

\theoremstyle{plain}
\newtheorem{lemma}[definition]{Lemma}
\newtheorem{proposition}[definition]{Proposition}
\newtheorem{theorem}[definition]{Theorem}
\newtheorem{corollary}[definition]{Corollary}
\newtheorem*{nntheorem}{Theorem}

\newenvironment{prf}{\begin{proof}[\textbf{\upshape Proof.}]}{\end{proof}}

\newenvironment{introduction}
	{\begin{center}{\textbf{\large Introduction}}\end{center}
	\markboth{Introduction}{Introduction}
}
	{}
\begin{document}
\sloppy
\title{\textbf{A Generalisation of Steenrod's Approximation Theorem}}
\renewcommand{\shortTitle}{Steenrod's Approximation Theorem}
\author{Christoph Wockel}
\date{}
\maketitle
\thispagestyle{empty}

\begin{abstract}
In this paper we aim for a generalisation of the Steenrod
Approximation Theorem from \cite[Section 6.7]{steenrod51}, concerning
a smoothing procedure for sections in smooth locally trivial
bundles. The generalisation is that we consider locally trivial smooth
bundles with a possibly infinite-dimensional typical fibre. The main
result states that a continuous section in a smooth locally trivial
bundles can always be smoothed out in a very controlled way (in terms
of the graph topology on spaces of continuous functions), preserving
the section on regions where it is already smooth.\\[\baselineskip]
\end{abstract}

\begin{introduction}
This paper generalises a result of \textsc{Steenrod} on a very nice
smoothing procedure for sections in locally trivial smooth bundles. It
puts together ideas from \cite[Section 6.7]{steenrod51}, \cite[Chapter
2]{hirsch76} and \cite[Section A.3]{neeb03} and tries to produce a
theorem of maximal generality out of them.

\begin{nntheorem}[Generalised Steenrod Approximation Theorem] Let $M$
be a finite-dimensional connected manifold with corners, $\pi \from
E\to M$ be a locally trivial smooth bundle with a locally convex
manifold $N$ as typical fibre and $\sigma \from M\to E$ be a
continuous section. If $L\se M$ is closed and $U\se M$ is open such
that $\sigma$ is smooth on a neighbourhood of $L\setminus U$, then for
each open neighbourhood $O$ of $\sigma (M)$ in $E$, there exists a
section $\tau \from M\to O$ which is smooth on a neighbourhood of $L$
and equals $\sigma$ on $M\setminus U$. Furthermore, there exists a
homotopy $F\from [0,1]\times M\to O$ between $\sigma$ and $\tau$ such
that each $F(t,\cdot)$ is a section of $\pi$ and $F(t,x)=\sigma
(x)=\tau (x)$ if $(t,x)\in [0,1]\times (M\setminus U)$.
\end{nntheorem}

This theorem is of maximal generality in the sense that the proof
depends heavily on the local compactness of $M$ an the local convexity
of $N$. Also, the topology kept in mind is the graph topology on
spaces of continuous function, which is rather fine, e.g., in
comparison to the compact-open topology. Thus there seems to be no
result of greater generality (e.g., for arbitrary base-spaces and
arbitrary fibres) which can be shown with the same method of proof.

The paper is organised as follows. The first definitions and remarks
introduce the setting of calculus on locally convex vector spaces and
manifolds (with and without corners) modelled on such spaces. We then
recall some basic constructions on the smoothing procedure for
continuous functions with values in locally convex spaces, which we
shall need in the proof of the main theorem. After having proved the
main theorem, we formulate some immediate consequences of it,
concerning the relation of smooth and continuous sections and
functions. 

Eventually, we arrive at the analogue result from
\cite{krieglMichor2002} in the locally convex setting, stating that
smooth and continuous homotopies into locally convex manifolds
agree. However, our method of proof is different from the one used in
\cite{krieglMichor2002} since it uses heavily the existence of charts
onto convex subsets, which are not available for a convenient manifold
in general. This result is quite interesting, because it has nice
applications in bundle theory \cite{equiv}.
\end{introduction}

\begin{definition}
\label{def:differentialCalculus}(cf. \cite{hamilton82},
\cite{milnor84} and \cite{gloecknerneeb}) Let $X$ and $Y$ be locally
convex spaces and $U\se X$ be open. Then $f\from U\to Y$ is called
\textit{continuously differentiable} or \textit{$C^{1}$} if it is
continuous, for each $v\in X$ the differential quotient
\[
df (x).v:=\lim_{h\to 0}\frac{1}{h}(f (x+hv)-f (x))
\]
exists and the map $df\from U\times X\to Y$ is continuous. For $n>1$
we, recursively define 
\[
d^{n}f(x).(v_{1},\ldots,v_{n}):=\lim_{h\to 0}\frac{1}{h}
\left(d^{n-1}f(x+h).(v_{1},\ldots,v_{n-1})-
d^{n-1}f(x).(v_{1},\ldots,v_{n})\right)
\]
and say that $f$ is \textit{$C^{n}$} if $d^{k}f:U\times X^{k}\to Y$
exists for all $k=1,\ldots,n$ and is continuous. We say that $f$ is
$C^{\infty}$ or \textit{smooth} if it is $C^{n}$ for all $n\in\N$.

From this definition, the notion of a \textit{locally convex manifold}
is clear, i.e., a Hausdorff space such that each point has
neighbourhood that is homeomorphic to an open subset of some locally
convex space such that the corresponding coordinate changes are
smooth. Together with such a fixed differentiable structure on $M$, we
speak of $M$ as a locally convex manifold.
\end{definition}

\begin{remark}
In order to relate our results to other frequently used concepts of
differential calculus on infinite-dimensional vector spaces and
infinite-dimensional manifolds, we shortly line out the relation to
our setting (cf. \cite{keller74} for a more exhaustive comparison,
where smooth maps in our setting are called $C^{\infty}_{c}$-maps). In
the case of Banach-spaces $X$ and $Y$, a map is called \emph{Fr\'echet
differentiable} if it is differentiable in the sense of Definition
\ref{def:differentialCalculus} and the differential $x\mapsto df(x)$
is a continuous map into the space of bounded linear operators
$B(X,Y)$, endowed with the norm topology. Thus, Fr\'echet
differentiable (resp. smooth) maps are differentiable (resp. smooth) in our
setting.

Next, we recall the basic definitions of the convenient
calculus from \cite{krieglmichor97}. Let $X$ and $Y$ be arbitrary locally
convex spaces.  A curve $f:\R\to X$ is called smooth if it is smooth
in the sense of Definition \ref{def:differentialCalculus}.  Then the
$c^{\infty}$-topology on $X$ is the final topology induced from all
smooth curves $f\in C^{\infty} (\R,X)$. If $X$ is a Fr\'echet space,
then the $c^{\infty}$-topology is again a locally convex vector
topology which coincides with the original topology \cite[Theorem
4.11]{krieglmichor97}.  If $U\se X$ is $c^{\infty}$-open, then $f:U\to
Y$ is said to be of class $C^{\infty}$ or smooth if
\[
f_{*}\left(C^{\infty} (\R,U) \right)\se C^{\infty} (\R,Y),
\]
i.e., if $f$ maps smooth curves to smooth curves. The chain rule
\cite[Proposition 1.15]{gloeckner02b} implies that each smooth map in
the sense of Definition \ref{def:differentialCalculus} is smooth in
the convenient sense. On the other hand, \cite[Theorem
12.8]{krieglmichor97} implies that on a Fr\'echet space a smooth map
in the convenient sense is smooth in the sense of Definition
\ref{def:differentialCalculus}. Hence for Fr\'echet spaces, this notion
coincides with the one from Definition \ref{def:differentialCalculus}.
\end{remark}

\begin{definition}
A \emph{$d$-dimensional manifold with corners} is a paracompact
Hausdorff space such that each point has a neighborhood that is
homeomorphic to an open subset of
\[
\R^{d}_{+}=\{(x_{1},\dotsc,x_{d})\in\R^{d}:x_{i}\geq 0\text{ for all
}i=1,\dotsc,d\}
\]
and such that the corresponding coordinate changes are smooth
(cf. \cite{lee03}). The crucial point here is the notion of smoothness
for non-open domains.  The usual notion is to define a map $f:A\se
\R^{n}\to \R^{m}$ to be smooth if for each $x\in A$, there exists a
neighborhood $U_{x}$ of $x$ which is open in $\R^{n}$, and a smooth
map $f_{x}:U_{x}\to \R^{m}$ such that $\left.f_{x}\right|_{A\cap
U_{x}}=\left.f\right|_{A\cap U_{x}}$.

A more general concept of manifolds with corners modeled on locally
convex spaces can be found in \cite{michor80}, \cite{smoothExt} and
\cite{gloecknerneeb}, along with the appropriate definitions of
differentiable or smooth functions. Basically, in this setting, a map
on a non-open domain with dense interior is defined to be smooth if it
is smooth on the interior and differentials extend continuously to the
boundary.
\end{definition}

\begin{remark}\label{topologyonmanifold}
We recall some basic facts from general topology. A topological space
$X$ is called \emph{paracompact} if each open cover has a locally
finite refinement. If $X$ is the union of countably many compact
subsets, then it is called \emph{$\sigma$-compact}, and if each open
cover has a countable subcover, it is called \emph{Lindel\"of}.

Now, let $M$ be a fi\-nite-\-di\-men\-sion\-al manifold with corners,
which is in particular locally compact and locally connected. For
these spaces, \cite[Theorems XI.7.2+3]{dugundji66} imply that $M$ is
paracompact if and only if each component is $\sigma$-compact,
equivalently, Lindel\"of.  Furthermore, since paracompact spaces are
normal, $M$ is normal in each of these cases.

One very important fact on $M$ is that it permits smooth partitions of
unity (c.f.\ \cite[Theorem 2.1]{hirsch76}). That means that for each
locally finite open cover $(V_{i})_{i\in I}$ we find smooth functions
$\lambda_{i}\from M\to [0,1]$ such that $\supp (\lambda_{i})\se V_{i}$
and $\sum_{i\in I}(\lambda_{i}(x))=1$.
\end{remark}

\begin{definition}
If $X$ is a Hausdorff space and $Y$ is a topological space, then
$C(X,Y)_{c.o.}$ is the space of continuous functions from $X$ to $Y$,
endowed with the compact-open topology (cf. \cite[Section
X.3.4]{bourbakiTop}). A basic open set in this topology is given by
$\lfloor C_{1},W_{1}\rfloor\cap \ldots\cap \lfloor C_{n},W_{n}\rfloor$
for $C_{1},\dots ,C_{n}\se X$ compact and $W_{1},\dots ,W_{n}\se Y$
open, where 
\[
\lfloor C,W\rfloor:=\{f\in C(X,Y): f(C)\se W\}.
\]
\end{definition}

\begin{remark}
If $Y$ happens to be a topological group, then this topology coincides
with the topology of compact convergence \cite[Theorem
X.3.4.2]{bourbakiTop} and thus $C(X,Y)_{c.o.}$ is a topological group
itself.  If $Y$ is a locally convex space, then $C(X,Y)_{c.o.}$ is
again a locally convex space space with respect to pointwise
operations.

If $Y$ is a locally compact
space, then the exponential law yields that the canonical map $\varphi
\from C(X,C(Y,Z))\to C(X\times Y,Z)$, $\varphi (f)(x,y)=f(x)(y)$ is a
homeomorphism \cite[Section X.3.4]{bourbakiTop}.

A finer topology on $C(X,Y)$ is the graph topology, which we term
$C(X,Y)_{\Gamma}$
(cf. \cite{naimpally66GraphTopologyForFunctionSpaces}). A basic open
set in this topology is given by $\Gamma_{U}:=\{f\in C(X,Y):\Gamma
(f)\se U\}$, where $\Gamma(f)$ denotes the graph of $f$ in $X\times Y$
and $U\se X\times Y$ is open.
\end{remark}

\begin{proposition}
If $M$ is a finite-dimensional $\sigma$-compact manifold with corners, then
for each locally convex space $Y$ the space $C^{\infty }(M,Y)$ is dense
in $C (M,Y)_{c.o.}$. If $f\in C (M,Y)$ has compact support and U is an open
neighbourhood of $\supp(f)$, then each neighbourhood of $f$ in $C (M,Y)$
contains a smooth function whose support is contained in $U$.
\end{proposition}
\begin{prf}
The proof of \cite[Theorem A.3.1]{neeb03} carries over without changes.
\end{prf}
\begin{corollary}\label{cor:vectorValuedFunctionsAreDense}
If $M$ is a finite-dimensional $\sigma$-compact manifold with corners and $V$ 
is an open subset of the locally convex space $Y$, then $C^{\infty} (M,V)$ is
dense in $C (M,V)_{c.o.}$.
\end{corollary}

\begin{lemma}\label{teclemma6}
Let $M$ be a finite-dimensional $\sigma$-compact manifold with
corners, $Y$ be a locally convex space, $W\se Y$ be open and convex
and $f:M\to W$ be continuous. If $L\se M$ is closed and $U\se M$
is open such that $f$ is smooth on a neighbourhood of $L\setminus U$,
then each neighbourhood of $f$ in $C (M,W)_{c.o.}$ contains a continuous
map $g:M\to W$, which is smooth on a neighbourhood of $L$ and which equals $f$
on $M\setminus U$.
\end{lemma}

\begin{prf}
(cf. \cite[Theorem 2.5]{hirsch76}) Let $A\se M$ be an open set
containing $L\setminus U$ such that $f\big|_{A}$ is smooth.  Then
$L\setminus A\se U$ is closed in $M$, and, since $M$ is normal
(cf. Remark \ref{topologyonmanifold}), there exists $V\se U$ open with
$L\setminus A\se V\se \ol{V}\se U$. Then $\{U,M\setminus\ol{V}\}$ is
an open cover of $M$, and there exists a smooth partition of unity
$\{\lambda_{1},\lambda_{2}\}$ subordinated to this cover. Then
\[
G_{f}:C(M,W)_{c.o.}\to C(M,W)_{c.o.},\;\; 
G_{f}(\gamma)(x)=\lambda_{1}(x)\gamma(x) 
+\lambda_{2}(x)f(x)
\]
is continuous since $\gamma \mapsto \lambda_{1}\gamma$
and $\lambda_{1}\gamma \mapsto \lambda_{1}\gamma+\lambda_{2}f$ are continuous.

If $\gamma$ is smooth on $A\cup V$ then so is $G_{f}(\gamma)$, because
$\lambda_{1}$ and $\lambda_{2}$ are smooth, $f$ is smooth on $A$ and
$\lambda_{2}\big|_{V}\equiv 0$. Note that $L\se A\cup (L\setminus
A)\se A\cup V$, so that $A\cup V$ is an open neighbourhood of
$L$. Furthermore, we have $G_f(\gamma)=\gamma$ on $V$ and $G_f(\gamma
)=f$ on $M\setminus U $. Since $G_{f}(f)=f$, there is for each open
neighbourhood $O$ of $f$ an open neighbourhood $O'$ of $f$ such that
$G_{f}(O')\se O$. By Corollary \ref{cor:vectorValuedFunctionsAreDense}
there is a smooth function $h\in O'$ such that $g:=G_{f}(h)$ has the
desired properties.
\end{prf}

\begin{lemma}\label{cor3}
Let $M$ be a finite-dimensional $\sigma$-compact manifold with
corners, $N$ be a smooth manifold, modelled on a locally convex space,
$W\se N$ be diffeomorphic to an open convex subset of the modelling
space of $N$ and $f:M\to W$ be continuous. If $L\se M$ is closed and
$U\se M$ is open such that $f$ is smooth on a neighbourhood of
$L\setminus U$, then each neighbourhood of $f$ in $C (M,W)_{c.o.}$
contains a map $g\from M\to W$ which is smooth on a neighbourhood of
$L$ and which equals $f$ on $M\setminus U$.
\end{lemma}

\begin{prf}
Let $\varphi :W\to \varphi (W)$ be the postulated diffeomorphism. If
$\lfloor C_{1},V_{1}\rfloor \cap \ldots\cap \lfloor
C_{n},V_{n}\rfloor$ is an open neighbourhood of $f\in C(M,W)_{c.o.}$,
then $ \lfloor C_{1},\varphi (V_{1})\rfloor \cap \ldots\cap \lfloor
C_{n}, \varphi (V_{n})\rfloor$ is an open neighbourhood of $\varphi
\circ f$ in $C(M,\varphi (W))_{c.o.}$. We apply Lemma \ref{teclemma6}
to this open neighbourhood to obtain a map $h$. Then
$g:=\varphi^{-1}\circ h$ has the desired properties.
\end{prf}

\begin{theorem}[Generalised Steenrod Approximation Theorem]
\label{thm:approximationTheorem} Let $M$ be a finite-dimensional
connected manifold with corners, $\pi \from E\to M$ be a locally
trivial smooth bundle with a locally convex manifold $N$ as typical
fibre and $\sigma \from M\to E$ be a continuous section. If $L\se M$
is closed and $U\se M$ is open such that $\sigma$ is smooth on a
neighbourhood of $L\setminus U$, then for each open neighbourhood $O$
of $\sigma (M)$ in $E$, there exists a section $\tau \from M\to O$
which is smooth on a neighbourhood of $L$ and equals $\sigma$ on
$M\setminus U$. Furthermore, there exists a homotopy $F\from
[0,1]\times M\to O$ between $\sigma$ and $\tau$ such that each
$F(t,\cdot )$ is a section of $\pi$ and $F(t,x)=\sigma (x)=\tau (x)$
if $(t,x)\in [0,1]\times (M\setminus U)$.
\end{theorem}

\begin{prf}(cf. \cite[Section 6.7]{steenrod51}) We describe roughly how
the proof is going to work. After choosing an appropriate cover
$(V_{i})_{i\in \N}$ of $M$ in the beginning, we shall inductively
construct sections $\tau_{i}$ of $\pi$, that become smooth on
increasing subsets of $M$. To avoid convergence considerations, we
construct $\tau$ stepwise from the $\tau_{i}$ in the end.

We claim that there exist locally finite open covers $(V_{i})_{i\in
\N_{0}}$, $(V'_{i})_{i\in N_{0}}$ of $M$, $(W_{i})_{i\in \N_{0}}$ of
$\sigma (M)$ and a collection $(Z_{i})_{i\in \N_{0}}$ of open subsets
of $N$ which are diffeomorphic to convex open subsets of the modelling
space of $N$, such that we have
\begin{itemize}
\item $\ol{V'_{i}}$ and $\ol{V_{i}}$ are compact
\item $\ol{V'_{i}}\se V_{i}$ and $W_{i}\se O$
\item $\sigma (\ol{V_{i}})\se W_{i}$ (which is equivalent to $
       \ol{V_{i}}\se \pi(W_{i})$ and implies $
           V_{i}\se \pi (W_{i})$)
\item the restricted bundle $\left.\pi\right|_{\pi (W_{i})}$ is
trivial and there exist smooth trivialisations $\varphi_{i}\from
\pi^{-1}(\pi (W_{i}))\to \pi (W_{i})\times Y$ such that
$\varphi_{i}(\pi^{-1}(V_{i}))=V_{i}\times Z_{i}$
\end{itemize}
for each $\in \N_{0}$. First, we recall the properties of the topology
on $M$ from Remark \ref{topologyonmanifold}. Now, let $(S_{t})_{t\in
T}$ be a trivialising cover of $M$ and \mbox{$\varphi_{t}\from
\pi^{-1}(S_{t})\to S_{t}\times N$} be the corresponding local
trivialisations. That means, each $\varphi_{t}$ is a diffeomorphism
satisfying $\pr_{1}(\varphi_{t}(x))=\pi (x)$ for all $x\in
\pi^{-1}(S_{t})$. Then each $x\in M$ is in $S_{t(x)}$ for some map
$M\ni x\mapsto t(x)\in T$. Furthermore, there exist open
neighbourhoods $U_{x}\se S_{t(x)}$ of $x$ and $Z_{x}\se N$ of
$\pr_{2}\big(\varphi_{t(x)}(\sigma (x))\big)\in N$ such that $Z_{x}$
is diffeomorphic to an open convex subset of the modelling space of
$N$ and $\varphi^{-1}_{(t(x))}(U_{x}\times Z_{x})\se O$. Since $M$ is
normal, each $x\in M$ has a relatively compact open neighbourhood
$V_{x}$ such that $\ol{V_{x}}\se U_{x}$ and
$\pr_{2}\big(\varphi_{t(x)}(\sigma (\ol{V_{x}}))\big)\se
Z_{x}$. Furthermore, let $V'_{x}$ be an open neighbourhood of $x$ such
that $\ol{V'_{x}}\se V_{x}$. As $M$ is paracompact, $(V_{x})_{x\in M}$
has a locally finite refinement $(V_{j})_{j\in J}$. Since each
$\ol{V_{j}}$ is covered by finitely many $V'_{x_{j,1}},\dots
,V'_{x_{j,k_{j}}}$, we deduce that
\[
(V_{j}\cap V'_{x_{j,1}},\dots,V_{j}\cap V'_{x_{j,k_{j}}})_{j\in J}
\]
is also a locally finite open cover of $M$. By re-defining the index
set we thus get two locally finite open covers $(V_{i})_{i\in I}$ and
$(V'_{i})_{i\in I}$ such that $\ol{V'_{i}}$ and $\ol{V_{i}}$ are
compact and we have $\ol{V'_{i}}\se V_{i}$ for each $i\in I$. In
addition, we may assume that $I=\N^{+}$ for $M$ is Lindel\"of.

Since each $\ol{V_{i}}$ is contained in some $\ol{V_{x(i)}}$ of $M$
and the values of $\pr_{2}\op{\circ}\varphi_{t(x(i))}\op{\circ}\sigma$
on $\ol{V_{x(i)}}$ are contained in $Z_{x(i)}$, we get local
trivialisations
$\varphi_{i}:=\left.\varphi_{t(x(i))}\right|_{\pi^{-1}(U_{x(i)})}$ and
open subsets $Z_{i}:=Z_{x(i)}$ of $N$ and
$W_{i}:=\varphi^{-1}_{t(x(i))}(U_{x(i)}\times Z_{i})$ satisfying all
requirements.

We set $V_{0}:=\emptyset$ and $V'_{0}:=\emptyset$, and observe that
$(V_{i})_{i\in \N_0}$ and $(V'_{i})_{i\in\N_{0}}$ are locally finite
covers by their construction. Furthermore, we assume that $\sigma$ is
smooth on the open neighbourhood $A$ of $L\setminus U$ and that $A'$
is another open neighbourhood of $L\setminus U$ with $\ol{A'}\se
A$. Define
\[
L_{i}:=L\cap \ol{V'_{i}}\setminus (V'_{0}\cup \ldots\cup V'_{i-1})
\]
which is closed and contained in $V_{i}$.  Since $L\setminus A'\se U$
we have $L_{i}\setminus A'\se V_{i}\cap U$ and there exist open
subsets $U_{i}\se V_{i}\cap U$ such that $L_{i}\setminus A'\se
U_{i}\se\ol{U_{i}}\se V_{i}\cap U$.  We claim that there exist
continuous sections $\tau_{i}\in C(M,E)$, $i\in \N_{0}$, satisfying
\begin{enumerate}
\renewcommand{\labelenumi}{(\theenumi)}
\renewcommand{\theenumi}{\alph{enumi}}
\item $\tau _{i}=\tau _{i-1}$ on $M\setminus \ol{U_{i}}$ for all
      $i\in \N^{+}$,\label{eqn:condA}
\item $\tau_{i}(M)\se O$ and $\tau _{i}(\ol{V_{j}})\se W_{j}$ for all 
      $i,j\in \N_{0}$,
      \label{eqn:condB}
\item $\tau _{i}$ is smooth on a neighbourhood of 
      $L_{0}\cup \ldots\cup L_{i}\cup\ol{A'}$ for all $i\in\N_{0}$ and
      \label{eqn:condC}
\item for each $i\in \N^{+}$ there exists a homotopy 
      $F_{i-1}\from[0,1]\times M\to O$ such that each $F(t,\cdot )$ is a 
      section of 
      $\pi$, $F_{i-1}(0,\cdot )=\tau_{i-1}$
      and $F_{i-1}(1,\cdot )=\tau_{i}$, which is
      constantly $\tau_{i-1}=\tau_{i}$ on $[0,1]\times (M\setminus \ol{U_{i}})$.
      \label{eqn:condD}
\end{enumerate}
Condition \eqref{eqn:condA} will ensure that we can construct $\tau$
stepwise from the $\tau_{i}$, and condition \eqref{eqn:condB} will
ensure that we can view $\left.\tau_{i-1}\right|_{V_{i}}$ as a
$Z_{i}$-valued function on $V_{i}$ and thus can apply Corollary
\ref{cor3} to $\left.\tau_{i-1}\right|_{V_{i}}$ in order to construct
$\tau_{i}$. Finally, condition \eqref{eqn:condC} will ensure the
asserted smoothness property of $\tau$, and condition
\eqref{eqn:condD} will enable us to construct the asserted homotopy.

For $i=0$ we set $\tau_{0}=\sigma$, which clearly satisfies conditions
\eqref{eqn:condA}-\eqref{eqn:condC}. Hence we assume that the $\tau
_{i}$ are defined for $i<a$. We consider the set
\[
Q:=\{\gamma\in C (V_{a},W_{a}):\gamma =\tau_{a-1}\text{ on }V_{a}\setminus 
\ol{U_{a}}\},
\]
which is a closed subspace of $C (V_{a},W_{a})_{c.o.}$. Then we have a
well-defined map
\[
G:Q\to C (M,W_{a}),\;\; G (\gamma) (x)=
\left\{\begin{array}{ll}
\gamma (x)&\text{if } x\in \ol{U_{a}}\\
\tau_{a-1} (x)&\text{if }x\in M\setminus \ol{U_{a}}.
\end{array} \right.
\]
Note that, by condition \eqref{eqn:condB}, we have
$\tau_{a-1}(V_{a})\se W_{a}$, whence $\left.\tau_{a-1}
\right|_{V_{a}}\in Q$. Furthermore, $\tau_{a-1}(M)\se O$ ensures
\begin{align}\label{eqn:neighbourhoodCondition}
G(\gamma)(M)\se O\text{\;\; if \;\;}
\gamma (V_{a})\se O.
\end{align}
Since $(V_{j})_{j\in
\N_{0}}$ is locally finite and $\ol{V_{j}}$ is compact, the set
$\{j\in \N_{0}:\ol{U_{a}}\cap \ol{V_{j}}\neq \emptyset\}$ is finite
and hence
\[
O'=\bigcap_{j\in \N_{0}}\lfloor \ol{U_{a}}\cap \ol{V_{j}},W_{a}\cap
W_{j}\rfloor
\]
is an open neighbourhood of $\left.\tau_{a-1} \right|_{V_{a}}$ in
$C(V_{a},W_{a})_{c.o.}$ by condition \eqref{eqn:condB}.

Since $\tau_{a-1}$ is a section, $\left.\tau_{a-1}\right|_{V_{a}}$ is
also a section of the restricted bundle
$\pi_{a}:=\left.\pi\right|_{V_{a}}$. For $\pi_{a}$ has the smooth
trivialisation $\varphi_{a}$, the space of sections of $\pi_{a}$ is
homeomorphic to $C(V_{a},N)$ by the homeomorphism $\sigma' \mapsto
H(\sigma '):=\pr_{2}\op{\circ }\varphi_{a} \op{\circ }\sigma '$ with
inverse given by $f\mapsto \varphi^{-1}_{a}\op{\circ}(\id\times
f)\op{\circ}\op{diag}$.  This shows in particular that $H(\sigma')$ is
smooth in a neighbourhood of $x\in V_{a}$ if and only if $\sigma'$ is
so.

We want to apply Lemma \ref{cor3} to
$g:=H(\left.\tau_{a-1}\right|_{V_{a}})\in C(V_{a},N)$ and claim for
this reason that $g$ takes values in some subset of $N$, diffeomorphic
to a convex neighbourhood of its modelling space. This in turn is
true, as $\varphi_{a}(W_{a})\se V_{a}\times Z_{a}$ and thus
$g=\pr_{2}\op{\circ}\varphi_{a}\op{\circ}\left.\tau_{a-1}\right|_{V_{a}}$
takes values in $Z_{a}$ by condition \eqref{eqn:condB}.

In order to construct $\tau_{i}$, we now apply Lemma \ref{cor3} to the
manifold with corners $V_{a}$, its closed subset $L_{a}':=(L\cap
\ol{V'_{a}})\cup (\ol{A'}\cap V_{a})\se V_{a}$, the open set $U_{a}\se
V_{a}$, $g\in C(V_{a},Z_{a})$ and the open neighbourhood $H(O')$ of
$g$.  Due to the construction, we have $L_{a}\setminus U_{a}\se A'\cap
V_{a}$ and, furthermore, $L\cap \ol{V'_{a}}\se L_{0}\cup \ldots \cup
L_{a}$. Hence we have
\[
    L_{a}'\setminus U_{a}
\se (L_{0}\cup \ldots\cup L_{a-1}\cup(L_{a}\setminus U_{a}))
    \cup (\ol{A'}\cap V_{a}\setminus U_{a})
\se L_{1}\cup \ldots\cup L_{a-1}\cup \ol{A'}
\]
so that by condition \eqref{eqn:condC}, $¸\left.\tau_{a-1}
\right|_{V_{a}}$ and, consequently, $g$ are smooth on a neighbourhood
of $L_{a}'\setminus U_{a}$. We thus obtain a map $h\in H(O')\se
C(V_{a},Z_{a})$ which is smooth on a neighbourhood of
$L_{a}'$. Furthermore, $H^{-1}(h)$ coincides with
$\left.\tau_{a-1}\right|_{V_{a}}$ on $V_{a}\setminus U_{a}\supseteq
V_{a}\setminus \ol{U_{a}}$, because there $h$ coincides with $g$ and
$h(x)=g(x)$ implies $H^{-1}(h)(x)=H^{-1}(g)(x)=\tau_{a-1}(x)$ for
$x\in V_{a}$. As a consequence, $H^{-1}(h)$ is contained in $O'\cap
Q$, and we set $\tau_{a}:=G (H^{-1}(h))$.

It remains to check that $\tau_{a}$ satisfies conditions
\eqref{eqn:condA}-\eqref{eqn:condD}.  Since $H^{-1}(h)$ coincides with
$\left.\tau_{a-1}\right|_{V_{a}}$ on $V_{a}\setminus \ol{U_{a}}$,
condition $\eqref{eqn:condA}$ is satisfied. From the construction we
know that $H^{-1}(h)(V_{a})\se W_{a}\se O$, which implies
$G(H^{-1}(h))(M)\se O$ by \eqref{eqn:neighbourhoodCondition}. In
addition $H^{-1}(h)\in O'$, which implies in turn
$H^{-1}(h)(\ol{U_{a}}\cap \ol{V_{j}})\se W_{j}$ and, furthermore,
$\tau_{a}(\ol{V_{j}})\se W_{j}$. Eventually, condition
\eqref{eqn:condB} is fulfilled. Furthermore, $\tau_{a}$ inherits the
smoothness properties from $\tau_{a-1}$ on $M\setminus \ol{U_{a}}$,
from $h$ on $V_{a}$ and since $L_{a}\se L\cap \ol{V'_{a}}$, condition
\eqref{eqn:condC} also holds. To construct $F_{a-1}$, we set
\[
F_{a-1}\from [0,1]\times M\to O,\quad (t,x)\mapsto 
\left\{\begin{array}{l@{\text{ if }}l}
\tau_{a-1}(x) & x\notin V_{a}\\
\varphi_{a}^{-1}(x,(1-t)\cdot g(x)+t\cdot h(x)) & x\in V_{a},
\end{array}\right.
\]
where the convex combination between $g(x)$ and $h(x)$ in $Z_{a}$ has
to be understood in local coordinates in the convex set which $Z_{a}$
is diffeomorphic to. Since $g$ equals $h$ on $V_{a}\setminus
\ol{U_{a}}$ and we have $\varphi_{a}^{-1}(x,g(x))=\tau_{a-1}(x)$ for
$x\in V_{a}$, this defines a continuous map which satisfies the
requirements of condition \eqref{eqn:condD}. This finishes the
construction of $\tau_{a}$ and thus the induction.

We next construct $\tau$. First we set $m (x):=\text{max}\{i:x\in
\ol{V_{i}}\}$ and $n (x):=\text{max}\{i:x\in V_{i}\}$. Then obviously
$n(x)\leq m(x)$ and each $x\in M$ has a neighbourhood on which
$\tau_{n(x)},\ldots,\tau_{m(x)}$ coincide since $\ol{U_{i}}\se V_{i}$
and $\tau_{i}=\tau_{i-1}$ on $M\setminus U_{i}$. Hence
$\tau(x):=\tau_{n(x)}(x)$ defines a continuous function on $M$. If
$x\in L$, then $x\in L_{0}\cup \ldots\cup L_{n(x)}$ and thus $\tau$ is
smooth on a neighbourhood of $x$. If $x\in M\setminus U$, then
$x\notin U_{1}\cup \ldots\cup U_{n(x)}$ and thus
$\tau(x)=\sigma(x)$. 

We finally construct the homotopy $F$. First observe that if $x\in M$
and $n> n(x)$, then $x\notin V_{n}\supseteq \ol{U_{n}}$ and
$F_{n}(t,x)=\tau_{n-1}(x)=\tau_{n-2}(x)=\dots =\tau_{n(x)}(x)=\tau
(x)$ by condition \eqref{eqn:condD}. We set
\[
F\from (0,1]\times M\to O,\quad (t,x)\mapsto F_{n-1}((n+1)\cdot
(1-nt),x) \text{ if }t\in[\frac{1}{n},\frac{1}{n+1}].
\]
This is well-defined and continuous since
$([\frac{1}{n},\frac{1}{n+1}]\times M)_{n\in \N^{+}}$ covers
$(0,1]\times M$, and for $n\geq 2$ and
\[
(t,x)\in  \left(\left[\frac{1}{n-1},\frac{1}{n}  \right]\times M \right)
     \cap \left(\left[\frac{1}{n},  \frac{1}{n+1}\right]\times M \right)
\]
we have $t=\frac{1}{n}$ and thus
$F(t,x)=F_{n-2}(1,x)=\tau_{n-1}(x)=F_{n-1}(0,x)$ by condition
\eqref{eqn:condD}. Furthermore, $\left.F\right|_{\{1\}\times
M}=F(0,\cdot )=\sigma$. We extend $F$ to
$[0,1]\times M$ by setting $F(0,x)=\tau(x)$. This is in fact a
continuous extension since each $x$ is contained in its open
neighbourhood $V_{m(x)}$ and $F(t,x')=\tau (x')$ for $(t,x')\in
[\frac{1}{m(x)},0]\times V_{m(x)}$ by the first observation of this
paragraph. Clearly, each $\left.F\right|_{\{t\}\times M}$ is a
section, because each $\left.F_{n}\right|_{\{t\}\times M}$ and $\tau$
are so. Furthermore, if $(t,x)\in [0,1]\times(M\setminus U)$, then
$x\notin U_{1}\cup \ldots\cup U_{n(x)}$ and thus $F_{0}(t,x)=\dots
=F_{n(x)}(t,x)=F(t,x)$.
\end{prf}

\begin{corollary}\label{cor:usualFormulation}
Let $M$ be a fi\-nite-\-di\-men\-sion\-al connected manifold with
corners, $N$ a be a locally convex manifold and $f\in C(M,N)$. If
$A\subseteq M$ is closed and $U\subseteq M$ is open such that $f$ is
smooth on a neighbourhood of $A\setminus U$, then each open
neighbourhood $O$ of $f$ in $C(M,N)_{\Gamma}$ contains a map $g$,
homotopic to $f$ in $O$, which is smooth on a neighbourhood of $A$ and
equals $f$ on $M\setminus U$. In particular, $C^{\infty}(M,N)$ is
dense in $C(M,N)_{\Gamma}$.

Furthermore, the same statement holds if we replace
the graph topology $C(X,Y)_{\Gamma}$ with the compact-open topology
$C(X,Y)_{c.o.}$.
\end{corollary}

\begin{prf}
We consider the globally trivial bundle $\pr_{1}\from M\times N\to M$.
Then the space of (continuous or smooth) mappings from $M$ to $N$ is
isomorphic to the space of (continuous or smooth) sections by
$f\mapsto \sigma_{f}$ with $\sigma_{f}(x)=(x,f(x))$. Then $\Gamma
(f)=\sigma_{f}(M)$ and the assertion follows directly from Theorem
\ref{thm:approximationTheorem} and the observation that the graph
topology is finer than the compact-open topology.
\end{prf}

\begin{proposition}\label{prop:homotopiesInSpacesOfSectionsAgree}
Let $M$ be a fi\-nite-\-di\-men\-sion\-al connected manifold with
corners and $\pi \from E\to M$ be a locally trivial smooth bundle with
a locally convex manifold $N$ as typical fibre. Then each continuous
section is homotopic to a smooth section. Furthermore, if there exists
a continuous homotopy between the smooth sections $\sigma$ and $\tau$,
then there exists a smooth homotopy between $\sigma $ and $\tau $.

Furthermore, we have that each base-point preserving continuous
section is homotopic, by a base-point preserving homotopy, to a
base-point preserving smooth section. Furthermore, if there exists a
continuous base-point preserving homotopy between the smooth sections
$\sigma$ and $\tau$, then there exists a smooth base-point preserving
homotopy between $\sigma $ and $\tau$.
\end{proposition}

\begin{prf}
The first assertion is already covered by Theorem
\ref{thm:approximationTheorem}. For the second assertion let $F'\from
[0,1]\times M\to E$ be a homotopy with $F'(0,\cdot)=\sigma $ and
$F'(1,\cdot )=\tau $.  Then we can construct a new homotopy $F''$
between $\sigma $ and $\tau $ which is smooth on a neighbourhood of
the closed subset $\{0,1\}\times M$ of the manifold with corners
$[0,1]\times M$. In fact, taking a smooth map $\gamma \from [0,1]\to
[0,1]$ with $\gamma ([0,\varepsilon ])=\{0\}$ and $\gamma
([1-\varepsilon ,1])=\{1\}$ for some $\varepsilon \in
(0,\frac{1}{2})$, $F''(t,x)=F'(\gamma (t),x)$ defines such a homotopy.
Now, $\sigma '\from [0,1]\times M\to [0,1]\times M\times E$,
$(t,x)\mapsto (t,x,F''(t,x))$ defines a section in the pull-back
bundle $\pr_{2}^{*}(E)$ of $E$ along the projection $\pr_{2}\from
[0,1]\times M\to M$. Furthermore, $\sigma '$ inherits the smoothness
properties of $F''$. Applying Theorem \ref{thm:approximationTheorem}
to the manifold with corners $[0,1]\times M$, the closed subset
$\{0,1\}\times M$, the open subset $(0,1)\times M$ and $\sigma '$
yields in the third component a smooth map $F\from [0,1]\times M\to E$
with $F(0,\cdot )=F''(0,\cdot )=F'(0,\cdot)=f$ and $F(1,\cdot
)=F''(1,\cdot )=F'(1,\cdot)=g$.

In the case of a base-point preserving section $\sigma$, we first
claim that there exists a base-point preserving homotopy to a section
which is constantly $f(x_{0})$ on a neighbourhood of the base-point
$x_{0}$ of $M$. In fact, if $x_{0}$ denotes the base-point of $M$,
then it has a neighbourhood $U$ such that there exists a local
trivialisation $\phi \from \pi^{-1}(U)\to U\times N$ and that
$\pr_{2}(\varphi (\sigma (U)))\se Z$ for some open subset $Z$ of $N$
which is diffeomorphic to the modelling space of $N$. We set
$f:=\left.\pr_{2}\op{\circ}\varphi \op{\circ}\sigma\right|_{U}$ and
take a smooth map $\lambda \from M\to [0,1]$ which is constantly $1$
on some neighbourhood of $x_{0}$ and with $\supp (\lambda )\se
U$. Then we define a homotopy
\[
G\from [0,1]\times M\to E,\quad
(t,x)\mapsto \left\{
\begin{array}{l@{\text{ if }}l}
\sigma (x) & x\notin U\\
\varphi^{-1}\big(x,\big(1-t\lambda (x)\big)\cdot f(x)+t\lambda(x)\cdot
f (x_{0})\big) & x\in U,
\end{array} \right.
\]
where the convex combination between $f(x)\in Z$ and $f(x_{0})\in Z$
has to be understood in local coordinates in the convex set which $Z$
is diffeomorphic to.  It is easily verified that $G$ is continuous and
has the desired properties, so we may assume that $\sigma$ is already
smooth on some neighbourhood of $x_{0}$. Again, interpreting $G$ as a
section in the pull-back bundle $\pr_{2}^{*}(E)$ as in the first part
of the proof and applying Theorem \ref{thm:approximationTheorem}
yields a homotopy which is constantly $\sigma (x_{0})$ on
$[0,1]\times\{x_{0}\}$.

Similarly, if $\sigma$ and $\tau$ are smooth and homotopic by a
continuous and base-point preserving homotopy $F'$, then the
construction of the first part of the proof yields a homotopy $F''$
which is yet base-point preserving. By a partition of unity argument
similar to construction on $G$, we may also assume that $F''$ is
continuously $\sigma (x_{0})=\tau (x_{0})$ on a neighbourhood of
$[0,1]\times \{x_{0}\}$. Once more, interpreting $F''$ as a section in the
pull-back bundle $\pr_{2}^{*}(E)$ as in the first part of the proof
and applying Theorem \ref{thm:approximationTheorem} yields a smooth
homotopy which coincides with $F''$ on $\{0,1\}\times M\cup
[0,1]\times \{x_{0}\}$ and thus meets all requirements.
\end{prf}

\begin{corollary}(cf. \cite{krieglMichor2002})
Let $M$ be a fi\-nite-\-di\-men\-sion\-al connected manifold with
corners and $N$ a be a locally convex manifold. Then each continuous
map is homotopic to a smooth map. Furthermore, two smooth maps $f$ and
$g$ are homotopic if and only if they are smoothly homotopic, i.e.,
there exists a smooth map $F\from [0,1]\times M\to N$ with $F(0,\cdot )=f$
and $F(1,\cdot )=g$.
\end{corollary}

\section*{Acknowledgements} The work on this paper was financially
supported by a doctoral scholarship from the Technische Universit\"at
Darmstadt. The author would like to express his thank to Karl-Hermann
Neeb and Helge Gl\"ockner for several very useful discussions on the
content of the paper.

\bibliography{mybib}
\vskip\baselineskip
\vskip\baselineskip
\vskip\baselineskip
\large
Christoph Wockel\\
Fachbereich Mathematik\\
Technische Universit\"at Darmstadt\\
Schlossgartenstrasse 7\\
D-64289 Darmstadt\\
Germany\\[\baselineskip]
\normalsize
\texttt{wockel@mathematik.tu-darmstadt.de}
\end{document}